\documentstyle[11pt]{article}
\input{epsf.sty}
\newtheorem {thm}{Theorem}[section]
\newtheorem {prop}[thm]{Proposition}
\newtheorem {lem}[thm]{Lemma}

\newcommand{\bea}{\begin{eqnarray}}
\newcommand{\ba}{\begin{array}}
\newcommand{\bean}{\begin{eqnarray*}}
\newcommand{\ea}{\end{array}}
\newcommand{\eea}{\end{eqnarray}}
\newcommand{\eean}{\end{eqnarray*}}
\newcommand{\be}{\begin{equation}}
\newcommand{\ee}{\end{equation}}

\newcommand{\lab}{\label}

\def\Cox{\hfill \Box}

\def\E{{\bf{E}}}
\def\P{{\bf{P}}}
\def\Nn{\hbox{I\kern-.2em\hbox{N}}}
\def\R{\hbox{I\kern-.2em\hbox{R}}}

\def\Z{{\bf{Z}}}

\def\|{\, | \, }

\def\0{\hat{0}}
\def\1{\hat{1}}

\def\th{\theta}
\def\Th{\Theta}

\def\N{N}

\begin{document}
\begin{center}
{\large \bf  No Directed Fractal Percolation in Zero Area} \\
\end{center}
\vspace{2ex}
\begin{center}
{\sc L. Chayes} \footnote{Department of Mathematics,
University of California, Los Angeles, California 90024.}, \, 
{\sc Robin Pemantle} \footnote{Research supported in part by 
National Science Foundation grant \# DMS 9300191, by a Sloan Foundation
Fellowship, and by a Presidential Faculty Fellowship.}$^,$\footnote{Department 
of Mathematics, University of Wisconsin, Madison, WI 53706 \, .}
\,  and \, 
{\sc Yuval Peres}\footnote{Research partially supported by NSF grant
\# DMS 9404391}$^,$\footnote{Mathematics Institute, The Hebrew University, 
Jerusalem, Israel
and Department of Statistics,  University of California, Berkeley, CA  \, .}
\end{center}
\noindent{\bf Abstract.}\, We consider the fractal percolation process on 
the unit square, with fixed decimation parameter $\N$  and level dependent 
retention parameters  $\{p_k\}$;
that is, for all $k \geq 1$, at the $k$th stage every retained square of side-length $\N^{1-k}$
is partitioned into $\N^2$ congruent subsquares, and each of these is retained
with probability $p_k$, independently of all others.
    We show that if $\prod_k p_k=0$ (i.e., if the area of the limiting set
vanishes a.s.) then a.s.\ the limiting set contains no directed crossings of 
the unit square
(a directed crossing is a path that crosses the unit square
from left to right, and moves only up, down and to the right).

\vfill




\section{Introduction}
Determining the most regular curves that are contained in the planar 
Brownian path $B[0,1]$ is a well-known open problem (see 
Duplantier et al (1988));
even the intuitively 
``obvious'' statement that $B[0,1]$ does not contain line segments is
not easy to prove (see Pemantle 1995) and it is not known whether
$B[0,1]$ contains  oriented paths or   directed paths.

Since this problem has so far been intractable,
in this paper we consider a ``mean field approximation''
to it, where the Brownian path $B[0,1]$ is
replaced by a random fractal with more statistical independence,
arising from ``fractal percolation''.

Fractal percolation is a process that generates a sequence $\{A_k\}$ of random
subsets of $[0,1]^2$. 
Let $\{p_k\}_{k \geq 1}$ be a sequence of numbers in 
$(0,1]$ and let $\N \geq 2$ be an
integer.  
The unit square, denoted $A_0$, is partitioned
into $\N^2$ congruent squares each of which is independently {\it retained} with probability
$p_1$ or {\it discarded} with probability $(1-p_1)$.  The closure of the retained, or
surviving, squares constitutes the set $A_1$. For $k>1$ The set $A_k \subset A_{k-1}$
 is generated by repeating this procedure, appropriately rescaled, on all the surviving 
squares of $A_{k-1}$, using the parameter $p_k$.  The limiting set  defined by 
$
A_\infty = \bigcap_k A_k
$
is the principal focus of study. The numbers of retained squares of different sizes form a
branching process in a varying environment; the offspring distribution at the $k\/$th
generation is Binomial$(\N^2, p_k)$ (see figure 1). 

\begin{figure}[hptp]
\epsfxsize=8cm
\hfill{\epsfbox{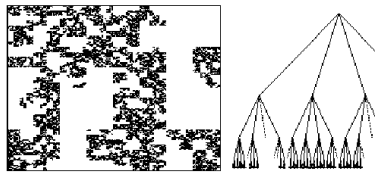}\hfill}
\caption{\bf $A_\infty$ and the corresponding branching process.}
\end{figure}

In the usual version of the
model, the probabilities $p_k$ are identically equal
to some  $p\in(0,1)$.
This is too crude for our purpose, since it yields random fractals of Hausdorff
dimension less than 2; to obtain a better approximation
to the Brownian path (which has dimension $2$), we will consider $p_k\to 1$;
this yields a limit set $A_\infty$ of Hausdorff dimension 2
(a.s.\ when $A_\infty$ is nonempty).

 For the case of constant retention probabilities $p_k=p$,
Chayes, Chayes and Durrett (1988)
showed that
$
\inf_k \P \{A_k \mbox{ \rm contains  a left-to-right crossing path } \}>0 \, ,
$
provided that $p$ is close enough to 1. 

For the standard percolation model in $\Z^d$, when the parameter
$p<1$ is large enough there is a.s.\ an infinite open  north-east
oriented path (see Durrett 1984). Recently, L. Chayes (1995) showed that for fractal
percolation with fixed $p<1$, there are a.s.\ 
no north-east oriented crossings of the unit square
in $A_\infty$, and furthermore, there are no {\em directed crossings}. 

\noindent{\bf Definition.}\, Let $\Gamma:[0,1] \mapsto[0,1]^2$ be a continuous path,
and write $\Gamma(t)=(\gamma_1(t),\gamma_2(t))$. The image $\Gamma[0,1]$ is called a 
(left-to-right) {\bf directed crossing} of the unit square if $\gamma_1(0)=0$,
the function $\gamma_1$ is nondecreasing, and  $\gamma_1(1)=1$.

Here we establish a stronger result:

\begin{thm} \label{thm:main}
Consider the 
fractal percolation process in $[0,1]^2$ with  decimation parameter
 $\N \geq 2$ and
retention parameters $\{p_k\}$.
If $\: \prod_{k=1}^{\infty} p_k = 0$, then a.s.\ 
$A_\infty$ does not contain a directed crossing of the unit square. 
\end{thm}

The condition  $\prod_k p_k = 0$ is clearly equivalent to 
the area of $A_\infty$ vanishing a.s., but this equivalence is
not used in the proof of the theorem.
The strategy of our proof is to exhibit many nearly horizontal ``contours'' 
(defined in the next section) that cannot be crossed
by a directed path in $A_\infty$. These contours force any directed crossing 
of the unit square 
in $A_\infty$ to be nearly horizontal; a simple counting argument precludes  
  nearly horizontal crossings.

The maximal height of the lowest contour at resolution $\N^{-r}$
is bounded by introducing the smallest vacancies first, and 
estimating the effect of larger vacancies by repeated application
of the following lemma:
\begin{lem} \lab{lem:expo}
Let $Y$ be an exponential random variable with parameter $\th$,
i.e., $\P[Y>y]= e^{-\th y}$ for $y>0$.
Then for any random variable $Z \geq 0$ which is 
independent of $Y$ and has finite mean, 
$$
\E\Big[\min\{Z,Y\}\Big] \leq \frac{2 \E Z}{2+ \th \E Z} \, .
$$
\end{lem}
\noindent{\bf Remark.} This inequality can also be used to relate
electrical conductance to expected maximum flow in a random network;
see Proposition \ref{prop:resist}.

\section{Proof of Theorem 1.1}
It will be convenient to extend the fractal percolation process to the
 half-strip
$[0,1] \times [0,\infty)$ by generating, in each of the squares  
$[0,1] \times [m-1,m]$ for $m \in \Nn$, an independent copy of the 
process as described in the introduction.

\noindent{\bf Definitions.} 
Let $k \geq 0$. A square $S$ of the form 
$[(j-1)N^{-k}, jN^{-k}] \times  [(\ell-1) N^{-k}, \ell N^{-k}]$,
 where $ 1 \leq j \leq N^k$ and $\ell \geq 1$, is called a
($k$th level) {\bf commensurate square}. \newline
In this section, it is convenient to construct $A_r$
in a different fashion than in the introduction. 
Let $\Xi_k$ denote a random subset of the collection of
level-$k$   commensurate squares, where each level-$k$ square
is in $\Xi_k$ with probability $1-p_k$, independently of all others.
Then $A_r$ is the union of all level-$r$  commensurate squares
that are not contained in any square from $\cup_{k=1}^r \Xi_k$. \newline 
 A  {\bf  contour} over the interval $[a,b] \subset [0,1]$
is a sequence of squares $S_1, \dots S_M$, each of which is in 
$\cup_{k=1}^\infty \Xi_k$, such that
\begin{itemize}
\item[(i)] the right edge of $S_i$ lies on the same line as 
the left edge of $S_{i+1}$ for $i = 1,\dots,M-1$ ;
\item[(ii)] the top of $S_{i+1}$ is at or above the same height as the bottom of 
$S_i$ for $i = 1,\dots,M-1$ ;
\smallskip
\item[(iii)]  the left edge of $S_1$ and the right edge of $S_M$ are on the lines $x = a$ and
$x = b$ respectively.
\end{itemize}

\begin{figure}[hptp]
\epsfxsize=6cm
\hfill{\epsfbox{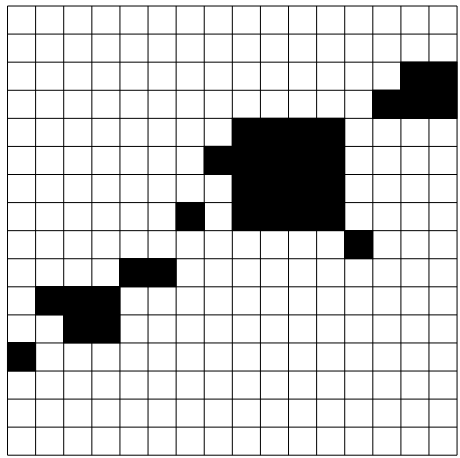}\hfill}
\caption{\bf  A contour over $[0,1]$.}
\end{figure}

The contour is called a {\em $(k,r)$-contour} 
if all the squares comprising it are
in $\cup_{j=k}^r \Xi_j$.
The bottom edges of the squares in a  
contour form the graph of a  step function;
the maximum of this function on the interval $[a,b] \subset [0,1]$ is called 
the {\bf height}
of the  contour over $[a,b]$.
\begin{prop} \lab{prop:compare}
 Fix an integer $r \geq 2$. Let $\th_k =\log(1/p_k)$ and 
$\Th_k = \sum_{j=1}^k \th_j$  for all $k$. 
 For $k \leq r$,
denote by $H_k$ the minimum, among all $(k,r)$-contours over $[0,\N^{-k}]$,
 of their height over that interval. Then 
\be \lab{eq:conb}
\E[H_k]  \leq \frac{2\N^{-k}}{\Th_r-\Th_{k-1}} \; \:
\mbox{ \rm  for all } k \leq r \,.
\end{equation}
 Moreover, the minimum, among all $(k,r)$-contours over $[0,1]$,
 of their height over $[0,1]$, has mean bounded above by 
$\frac{2}{\Th_r-\Th_{k-1}}$.
\end{prop}

\noindent{\bf Remark:} The contour heights $H_k$ that appear in the proposition
 also depend on $r$, but since $r$ is fixed throughout the proof of the 
proposition, it is suppressed from the notation.
 
To motivate the proof, consider momentarily
the cases where $\N=2$ and $ k$ is either $r$ or $r-1$. 
 Start by examining the process on $[0,2^{-r}] \times [0,\infty)$ and 
 considering only vacancies of scale $2^{-r}$.  The probability of observing a
``contour'' at height $\ell 2^{-r} $ is evidently $(1-p_r)(p_r)^{\ell}$, so
 $2^r H_r$ has expectation $p_r/(1-p_r) < \th_r^{-1}$.
  Now examine the situation on $[0,2\cdot2^{-r}]\times[0,\infty)$,
 illustrated in Figure 3. 

\begin{figure}[hptp]
\epsfxsize=6cm
\hfill{\epsfbox{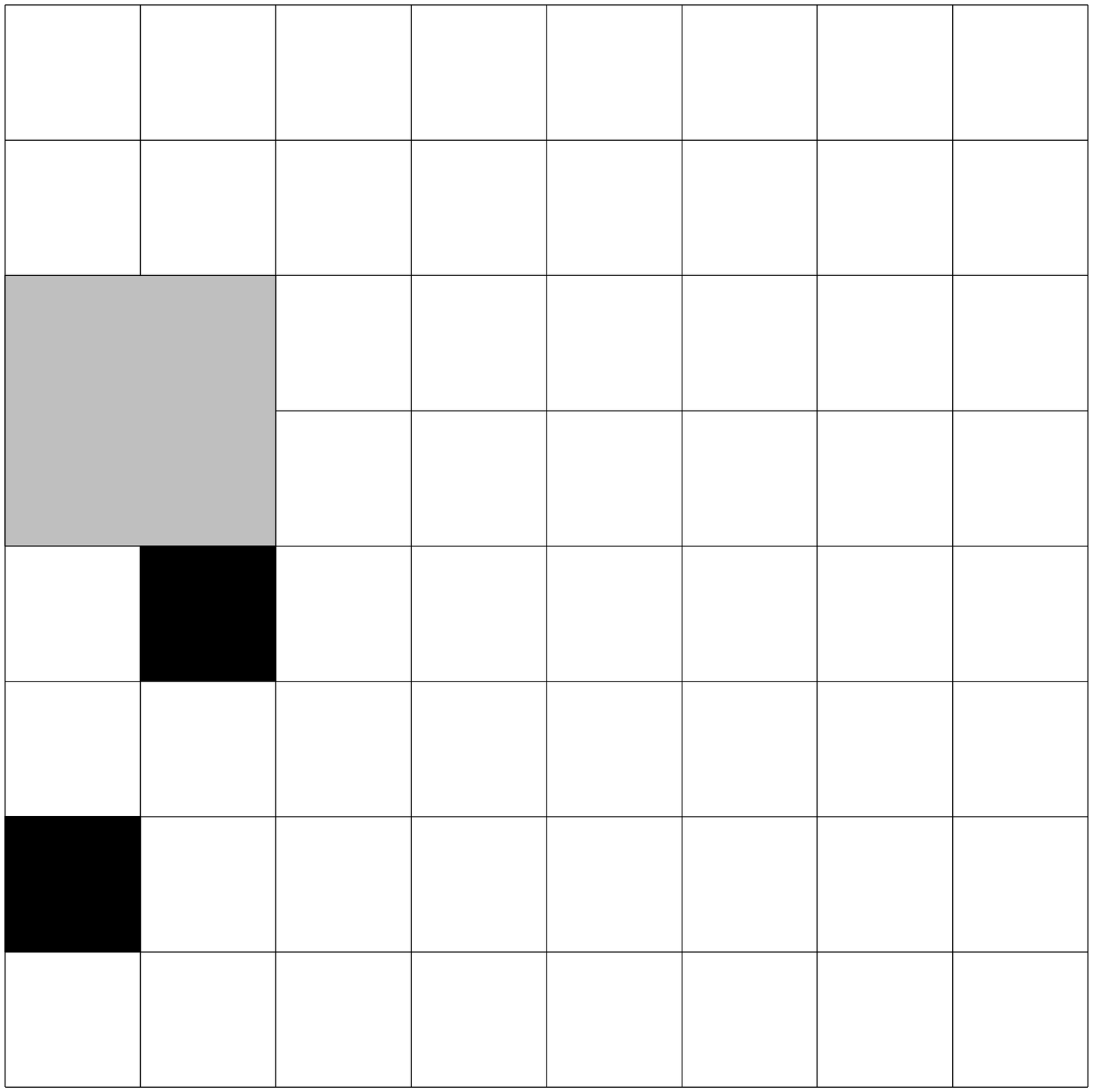}\hfill}
\caption{$H_{r-1} \leq 
\min \left\{ (H_{r} + H^\prime_{r}),\ 2^{1-r}  G_{r-1}  \right\}$ }
\end{figure}

Taking into account only the vacancies from $\Xi_r$, 
in the first column
we find the bottom edge of the lowest vacant square at height
$H_r$.  Assuming that $H_r >0$, we
may start the second column at height $H_r-2^{-r}$,
 and   go up a random height $H^\prime_{r}$ until reaching the next 
vacancy  from $\Xi_r$, where  $H_r$ and $H^\prime_{r}$ are i.i.d.  
However, the lowest square in $\Xi_{r-1}$ occurs
after $G_{r-1}$ squares of of side $2^{1-r}$, where
 $\, \P[G_{r-1} \geq \ell]=p_{r-1}^\ell$ for all $\ell \geq 0$. 
Recalling that $H_{r-1}$ is the 
height of the lowest $(r-1,r)$-contour  over   $ [0,2^{-r-1}] $, we have
\be \lab{eq:h2}
H_{r-1} \leq \min \left\{ (H_{r} + H^\prime_{r}),\ 2^{1-r} G_{r-1} 
 \right\}.
\ee
This inequality can be generalized to arbitrary decimation parameters
 $\N \geq 2$ 
and to coarser scales, but some care is needed, as  for levels $k<r$,
 contour heights need not be  integer multiples of the current scale $\N^{-k}$.
It is easily checked that $G_{r-1}$ is stochastically dominated by an 
exponential variable with parameter $\th_{r-1}$. Thus from (\ref{eq:h2}) and
Lemma \ref{lem:expo} we get
$$
\E[H_{r-1}] \leq \frac{4 \E [H_r]}{2+2^{r-1} \th_{r-1} \cdot 2 \E [H_r]}
\leq \frac{2 \cdot 2^{1-r}}{\th_r+\th_{r-1}} \,,
$$
a special case of (\ref{eq:conb}).

\noindent{\sc Proof of Proposition \ref{prop:compare}:}\,
Recall that $r \geq 2$ is fixed.
For any $v \geq 0$ and $k \leq r$, consider the class of 
 $(k,r)$-contours over $[0,\N^{-k}]$  such that the top of their 
leftmost square is at height  $\geq v$; 
denote by $H_k(v)$ the minimum, among these contours,
of their height above the line $y=v$ over the interval $[0,\N^{-k}]$.
We will prove that for all  $k \leq r$,
\be   \lab{eq:stoch} 
\E [H_k(v)] \leq  \frac{2\N^{-k}}{\Th_r-\Th_{k-1}}
\, \mbox{ \rm for all } \,
 v \geq 0 \, .
\end{equation}
The  proof proceeds inductively starting from $k=r$, and considering
the smallest (level-$r$) vacancies first.
Let $G_k(v)$ denote the number of of consecutive level-$k$
 squares on the leftmost column before the lowest square in $\Xi_k$,
counting upward from the level-$k$ square that intersects the line $y=v$.
(If there are two such squares, start from the higher one).
 Let $Y_k(v)$ be an exponential variable with parameter $\th_k$.
For integer $\ell \geq 0$, clearly
\be \lab{eq:geo}
\P[G_k(v) \geq  \ell] = p_k^{\ell}
 = e^{-\th_k \ell}= \P[\lfloor Y_k(v) \rfloor \geq \ell] \, , 
\end{equation}
so by enlarging the underlying probability space 
we can couple $G_k(v)$ with $Y_k(v)$ and assume that 
$G_k(v)=\lfloor Y_k(v) \rfloor$. Thus 
\be \lab{eq:eff}
H_k(v) \leq \N^{-k}G_k(v) \leq  \N^{-k}Y_k(v) 
 \mbox{ \rm for } 1 \leq k \leq r \,.
\end{equation}
In particular, this gives (\ref{eq:stoch}) for $k=r$.
For $k<r$ we can concatenate $\N \: \,$ 
 $(k+1,r)$-contours over intervals of length
$\N^{-k-1}$ to form a $(k,r)$-contour over $[0,\N^{-k}]$. 
This leads to the inequality
\be \lab{eq:sum}
H_k(v) \leq \sum_{i=1}^\N H_{k+1}^{(i)}(v_i) 
 \mbox{ \rm for }   k =0, \ldots, r-1   \, ,
\end{equation}
where $H_{k+1}^{(i)}(\cdot)$ measures contour heights over the interval 
$[\frac{i-1}{\N^{k+1}},\frac{i}{\N^{k+1}}]$, and the $v_i$ are defined by 
 $\, v_1=v$ and $v_{i+1}=v_i+H_{k+1}^{(i)}(v_i)$ 
(we suppress the dependence of $v_i$ on $v, k$ and $r$ from the notation).
Denote $Z_k(v)=\sum_{i=1}^\N H_{k+1}^{(i)}(v_i)$.
 Combining (\ref{eq:sum}) and the effect (\ref{eq:eff}) of level-$k$ vacancies,
we infer that
\be \lab{eq:min} 
H_k(v) \leq  \min\{Z_k(v), \N^{-k} Y_k(v)\}  \; 
\mbox{ \rm for }   k =1, \ldots, r-1 \, . 
\end{equation}
Note that  $H_{k+1}^{(i)}(\cdot)$ is determined by 
the collection of squares from $\cup_{\ell=k+1}^r \Xi_\ell$ in the strip 
$[\frac{i-1}{\N^{k+1}},\frac{i}{\N^{k+1}}] \times [0,\infty)$,
while $v_i$ is determined by the collection of squares from 
$\cup_{\ell=k+1}^r \Xi_\ell$ in  
the strip $[0, \frac{i-1}{\N^{k+1}}] \times [0,\infty)$; on the other hand,
 $G_k(v)$ only depends on $\Xi_k$, so we may take
$Z_k(v)$ and $Y_k(v)$ to be independent.

The induction hypothesis gives
$$ 
\E Z_k(v)= \sum_{i=1}^\N \E H_{k+1}^{(i)}(v_i)
 \leq  \frac{2\N^{-k}}{\Th_r-\Th_{k}} \, .
$$
By Lemma \ref{lem:expo} and the inequality (\ref{eq:min}),
$$
\E H_k(v) \leq \frac{2 \E[Z_k(v)]}{2 +\N^k \th_k \E[Z_k(v)]}
\leq  \frac{2 \N^{-k}}{\Th_r-\Th_k+ \th_k} \, .
$$
 This suffices to verify (\ref{eq:stoch}) by induction.

The last assertion of the proposition follows by concatenating $\N^k \: \,$
$(k,r)$-contours over intervals of length $\N^{-k}$.
$\Cox$

\noindent{\bf Proof of Theorem \ref{thm:main} completed:} \,
We are given that 
\be \lab{eq:given}
\Theta_r=\sum_{k=1}^r \log (1/p_k) \to \infty
\mbox{ \rm as } r \to \infty \, . 
\end{equation}
In particular
\be \lab{eq:partic}
1-p_k >\N^{-k/6} \mbox{ \rm for infinitely many } k\,.
\end{equation}
The probability that $A_\infty$ contains a specific point is at most
$4\prod_{k=1}^\infty p_k =0$, so we may restrict attention to directed
crossings that do not contain any point with both coordinates rational.
Such a crossing cannot pass through a contour. Fix $k<r$, and
consider $(k,r)$-contours over $[0,1]$ 
that lie above the line $y= \frac{i}{\N^k}$. Observe that if $S$
 is a  square commensurate at some level $\geq k$, and its bottom
is below the line $y= \frac{i+1}{\N^k}$, then $S$ must be contained in the
rectangle  $[0,1] \times [\frac{i}{\N^k},\frac{i+1}{\N^k}] $.
Proposition \ref{prop:compare} implies that the 
expected height  over $[0,1]$ above the line $y= \frac{i}{\N^k}$
of the lowest such $(k,r)$-contour is at most 
$\frac{2}{  (\Theta_{r}-\Theta_{k-1} )}$. By Markov's inequality,
the probability
that the rectangle $[0,1] \times [\frac{i}{\N^k},\frac{i+1}{\N^k}] $ 
does not contain a $(k,r)$-contour over $[0,1]$ is at most
$\frac{2\N^{k}}{ \Theta_{r} -\Theta_{k-1} } $.  This is
also an upper bound for the probability that the rectangle
$[0,1] \times [\frac{i-2}{\N^k},\frac{i-1}{\N^k}] $ does not contain
  a reflected
(``upside-down'') $(k,r)$-contour, provided that  $i\geq 2$. 
If the contour and the reflected contour just mentioned exist,
then a vertically aligned triple of  
 squares from $\Xi_k$ at any one of $\N^k$ locations 
(centered on row $i$ from the bottom) 
suffices to preclude the existence of a 
directed crossing that starts in the square
$[0, \frac{1}{\N^k}] \times [\frac{i-1}{\N^k},\frac{i}{\N^k}]$ 
and stays within $A_r$ (see Figure 4). For $i=1$, 
existence of the upper contour mentioned above and two vertically aligned
squares from $\Xi_k$ (in rows 1 and 2) 
 are enough to yield the same conclusion.

\begin{figure}[hptp]
\epsfxsize=6cm
\hfill{\epsfbox{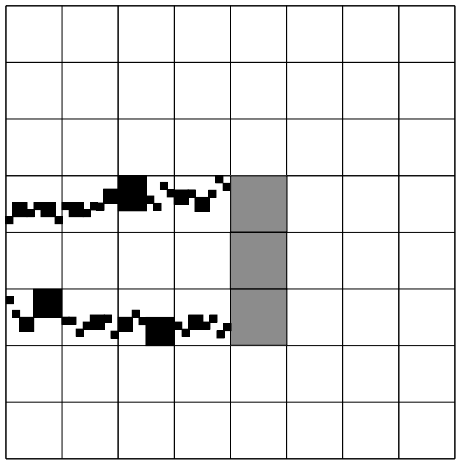}\hfill}
\caption{\rm A contour, a reflected contour, and three more vacancies,
  block directed  paths.} 
\end{figure}

Thus the probability that $A_\infty$ contains a directed 
crossing of the unit square is at most
$$
\lim_{r \to \infty}
\N^k\left(\frac{4 \N^{k}}{  \Theta_{r} -\Th_{k-1} }
 + (1-(1-p_k)^3)^{\N^k} \right)
\leq \N^k \exp(-\N^k (1-p_k)^3) 
$$
by (\ref{eq:given}). Finally, by (\ref{eq:partic}), for infinitely many $k$
 this probability is bounded by $\N^k \exp(-\N^{k/2}) $, so it must vanish.
$\Cox$
\section{Proof of Lemma 1.2}
For $z \geq 0$  we have 
$$
\E\Big[\min\{z,Y\}\Big] \, = \, \int_0^\infty \P[\min\{z,Y\}>y]  \, dy
=\int_0^z   e^{-\th y} \, dy \,  = \frac{1}{\th}(1- e^{-\th z}) \,.
$$
Hence for any nonnegative random variable $Z$, Jensen's inequality yields
\be \lab{eq:lap}
\E\Big[\min\{Z,Y\}\Big] = \frac{1}{\th} \E (1- e^{-\th Z}) 
\leq \frac{1}{\th} (1- e^{-\th \E Z}) \,.
\ee

Now rewrite the inequality
$$
2+x \geq \sum_{k=0}^\infty \frac{2-k}{k!} x^k =(2-x)e^{x}   
 \quad \forall x \geq 0
$$
in the equivalent form 
$$
1-e^{-x} \leq \frac{2x}{2+x} \,.
$$
Combining this with (\ref{eq:lap}) proves the lemma.
$\Cox$

\section{Concluding remarks}
\begin{enumerate}
\item An examination of the proof of Theorem \ref{thm:main}
has led to the following proposition concerning flow in a random network.
\begin{prop}[Lyons, Pemantle and Peres, 1996] \lab{prop:resist}
\, Consider a tree $\Gamma$ 
whose edges  are assigned independent exponential random
variables $\{Y(e)\}$ with means $\{m(e)\}$. Regard  $\{Y(e)\}$
as an assignment of capacities to the edges of $\Gamma$,
 and let $F$ denote the strength of the
maximum flow from the root to the boundary of $\Gamma$
 that is bounded by these capacities.
 On the other hand, let $\: {\cal C}$ be the 
effective conductance
from the root to the boundary for the edge conductances $m(e)$. 
Then $ \: {\cal C} \le \E[F] \le 2{\cal C}$.
\end{prop}
The maximum flow $F$ can also be represented as a first-passage 
percolation time in a certain planar dual of the tree;
 the lower bound
${\cal C} \le \E[F]$ (but not the upper bound)
extends to networks that are not trees.

\item We do not know if the ``zero area''
condition in Theorem \ref{thm:main} is 
necessary to prevent $A_\infty$ from containing
directed crossings. In other words, if $\prod_{k=1}^\infty p_k >0$,
does $A_\infty$
contain a directed crossing with positive probability?
\item 
 The Brownian path
is ``intersection-equivalent'' to the limiting set $A_\infty$
arising from fractal percolation with retention probabilities $p_k=k/(k+1)$,
see Peres (1996), but unfortunately intersection-equivalence
is too weak an equivalence relation to derive any rigorous conclusion
about Brownian motion from Theorem \ref{thm:main}. 
\end{enumerate}

\end{document}